\documentclass{amsart}
\usepackage{amssymb}
\usepackage{amscd}
\usepackage{verbatim}
\usepackage{epsfig}
\begin{document}

\newcommand\pa{\partial}
\newcommand\Ker{\operatorname{Ker}}
\newcommand\supp{\operatorname{supp}}
\newcommand\calC{{\mathcal C}}
\newcommand\Cinf{{\mathcal C}^{\infty}}
\newcommand\dist{{\mathcal C}^{-\infty}}
\newcommand\dCinf{\dot\Cinf}
\newcommand\ddist{\dot\dist}
\newcommand\ff{\operatorname{ff}}
\newcommand\mf{\operatorname{mf}}
\newcommand\bl{{\text b}}
\newcommand\scl{{\text{sc}}}
\newcommand\sci{{}^\scl}
\newcommand\Isc{I_\scl}
\newcommand\sfl{{\operatorname{s\Phi}}}
\newcommand\sfi{{}^\sfl}
\newcommand\WFsc{\operatorname{WF}_{\text{sc}}}
\newcommand\Osfh{\sfi\Omega^{\half}}
\newcommand\Isf{I_\sfl}
\newcommand\Vf{{\mathcal V}}
\newcommand\Vb{{\mathcal V}_{\bl}}
\newcommand\Vsc{{\mathcal V}_{\scl}}
\newcommand\Vsf{\Vf_\sfl}
\newcommand\sfT{{}^\sfl T}
\newcommand\scN{\sci N}
\newcommand\scT{\sci T}
\newcommand\yb{\bar y}
\newcommand\yt{\tilde y}
\newcommand\Yt{\tilde Y}
\newcommand\mut{\tilde \mu}
\newcommand\Lt{\tilde L}
\newcommand\RR{{\mathbb{R}}}
\newcommand\Lap{\Delta}
\setcounter{secnumdepth}{3}
\newtheorem{lemma}{Lemma}[section]
\newtheorem{prop}[lemma]{Proposition}
\newtheorem{thm}[lemma]{Theorem}
\newtheorem{cor}[lemma]{Corollary}
\newtheorem{result}[lemma]{Result}
\newtheorem*{thm*}{Theorem}
\newtheorem*{prop*}{Proposition}
\newtheorem*{conj*}{Conjecture}
\numberwithin{equation}{section}
\theoremstyle{remark}
\newtheorem{rem}[lemma]{Remark}
\theoremstyle{definition}
\newtheorem{Def}[lemma]{Definition}
\newtheorem*{Def*}{Definition}

\title{Intersecting Legendrians and blow-ups}
\author{Andrew Hassell}
\address{Centre for Mathematics and its Applications, Australian National
  University, Canberra ACT 0200 Australia}
\email{hassell@maths.anu.edu.au}
\author[Andr\'as Vasy]{Andr\'as Vasy}
\address{Department of Mathematics, Massachusetts Institute of Technology,
Cambridge MA 02139}
\email{andras@math.mit.edu}
\date{December 12, 2000}
\begin{abstract}
The purpose of this note is to describe the relationship between two
classes of Legendre distributions. These two classes are
distributions associated to an intersecting pair of Legendre submanifolds,
introduced
in \cite{Hassell:Plane} by analogy with intersecting Lagrangian distributions
of Melrose and Uhlmann \cite{Melrose-Uhlmann:Intersection}, and
Legendre distributions associated to a fibred scattering structure
introduced in \cite{Hassell-Vasy:Spectral}.
We prove a general result, and also give an example in two dimensions,
which shows explicitly the relation
between the two spaces in a simple setting.
\end{abstract}
\maketitle

\section{Introduction}

The purpose of the present note is to clarify the relationship between
two classes of Legendre distributions. The first class, that of
intersecting Legendrians associated to a pair of Legendre manifolds
which intersect cleanly, was defined by one of us in \cite{Hassell:Plane}
as an analog of the notion of intersecting Lagrangian distributions
\cite{Melrose-Uhlmann:Intersection} of Melrose and Uhlmann.
Just as the latter played an important role in the study of real
principal type operators, the former proved useful in geometric scattering
theory, both in describing the structure of the boundary value
of the resolvent (of a scattering Laplacian)
at the real axis \cite{Hassell-Vasy:Spectral}
and in the study of three-body scattering \cite{Hassell:Plane}.

The second class of Legendre distributions, that of
Legendrians associated to a fibred-scattering structure, was defined
in \cite{Hassell-Vasy:Spectral}. This extends the notion of Legendre
distributions to manifolds with corners that are equipped with
certain boundary fibrations, and it was used to analyze the structure of the
resolvent of scattering differential operators
near the the corners of the b-double space.

Our result is then that, given appropriate geometry and symbolic orders, the
class of intersecting Legendre distributions is a proper subset of
Legendre distributions associated to a fibred scattering structure --- see
Theorem~\ref{thm:2} in Section~\ref{sec:thm} for the precise statement. The
proper
inclusion corresponds to, roughly speaking, half of the possible terms
in a Taylor series expansion of a general fibred-scattering Legendre
distribution {\em not} being present in an intersecting Legendre
distribution.

The authors are grateful to Rafe Mazzeo and Richard Melrose for helpful
discussions. They acknowledge the support of the Australian Research Council
(A.\ H.) and the National Science Foundation through NSF grant \#DMS-99-70607
(A.\ V.).

\

Lagrangian distributions on a manifold without boundary, $X_0$, are
distributions on $X_0$ with very special singularities (in the
sense of lack of smoothness) associated to
Lagrangian submanifolds of the cotangent bundle of $X$. The simplest
examples are conormal distributions to an embedded submanifold $Z\subset X_0$;
these are distributions whose regularity is maintained under repeated
differentiation by vector fields tangent to $Z$; in particular they are
smooth away from $Z$.
Lagrangian distributions, and their generalizations, play a central
role in modern PDE theory, see e.g.\ \cite{Hor}.

If $X_0$ is not compact, one can study decay/growth properties of
distributions at `infinity' in addition to studying their smoothness
properties. Thus, the lack of rapid decay at infinity can be
considered a `singularity' and studied via microlocal analysis.
Since
one needs some structure at infinity, even to make sense of `rapid decay',
it is more natural to work on compact
manifolds with boundary (or corners) which arise by the compactification
of such $X_0$ and study singularities at the boundary.
On manifolds with boundary, $X$, one can either
introduce Legendre distributions from the symplectic (or really contact)
point of view, as traditionally done for Lagrangian distributions, or instead
simply write down such distributions and `work backwards'. Since
the former point of view, which is certainly `neater', has been discussed
in detail in \cite{RBMZw}, we follow the second approach. This should
quickly make it clear that many familiar functions fall in the class of
Legendre distributions.

Thus, let $X$ be a compact manifold with boundary, and $x$ a boundary defining
function. Legendre distributions on $X$ are functions which are smooth
in the interior of $X$ and have specific types of singularities at the
boundary, which we shall now describe. The simplest example is a function
of the form
\begin{equation}
u_1 = x^q e^{i\phi(x,y)/x} a(x,y), \text{ with } \phi, a \text{ smooth on } X.
\label{eq:Leg-1}\end{equation}
Here $(x,y)$ are coordinates on $X$, where $y = (y_1, \dots, y_{n-1})$
restrict to coordinates on $\pa X$. There is no loss of generality in
assuming that $\phi$ depends only on $y$. Another example is a
function
$$
x^q V(x, y/x),
$$
where $V(x,w)$ is smooth in $x$ and Schwartz in $w$. More generally, if $y
= (y', y'')$ is a splitting of the coordinates then a function of the form
\begin{equation}
u_2 = x^q V(x,y'/x,y'')
\end{equation}
with $V(x,w',y'')$ smooth and Schwartz in $w'$, is a Legendre
distribution. Distributions such as $u_1$ arise, for example, as plane
waves in scattering theory, while distributions like $u_2$ arise as
potentials in many-body scattering, as we explain below.
This already makes it clear that
we will need a combination of these types of distributions, e.g.\ to
understand perturbed plane waves, in a precise manner, in many-body scattering.

\

These two types of Legendre distributions can be given a uniform
description in terms of Legendre submanifolds. Let $U$ be an open subset of
$\pa X$, and consider the space $U \times \RR \times \RR^{n-1}$ with
coordinates $y, \tau, \mu$. The form $\chi = d\tau + \mu \cdot dy$ is a
contact form on $U \times \RR \times \RR^{n-1}$; that is, $\chi \wedge
(d\chi)^{n-1}$ never vanishes. (Here we work in local coordinates, but
this has an invariant geometric description in terms of the scattering
cotangent bundle; see \cite{RBMZw, Hassell-Vasy:Spectral} for a detailed
description, and the next section for a brief summary.)
Recall that a Legendre distribution is a
submanifold of maximal dimension (equal to $n-1$) on which the contact form
vanishes. The Legendre submanifold associated with $u_1$ is
$$
G_1 = \{ (y, \tau = -\phi(y), \mu = d_y \phi(y) ) \}.
$$
This submanifold determines $\phi$. For each $m \in \RR$, the class $\Isc^m(X,
G_1)$ is the class of functions $u_1$ of the form \eqref{eq:Leg-1} with $q
= m + n/4$ and $a \in \Cinf(X)$ arbitrary.

The second example, $u_2$, can be written in terms of the Fourier transform
in the $w$ variable as
\begin{equation}
u_2 = x^q \int e^{iy' \cdot \eta'/x} \hat V(x,\eta',y'') d\eta'.
\label{eq:Leg-2}
\end{equation}
The Legendre submanifold associated with $u_2$ is
$$
G_2 = \{ (y, \tau, \mu) \mid y' = 0, \tau = 0, \mu'' = 0 \}.
$$
If $\dim y' = k$ then the class $\Isc^m(X,G_2)$ is the class of functions of
the form \eqref{eq:Leg-2} with $q = m+n/4-k/2$ and with $\hat V$ Schwartz
in the second variable.

More generally, any Legendre submanifold $L$ of $U \times \RR \times \RR^{n-1}$
has a local parametrization, that is, a function $\phi(y,v)$, where $y \in
U$ and $v \in \RR^p$ such that locally
$$
L = \{ (y, \tau, \mu) \mid \exists \ (y,v) \text{ such that } \tau = -\phi(y,v), \mu = d_y \phi(y,v) \text
{ and } d_v \phi(y,v) = 0 \},
$$
and
\begin{equation}
\text{ for } 1 \leq i \leq p, \quad d \left( \frac{\pa \phi}{\pa v_i} \right)
\text{ are linearly independent. }
\label{nondeg}\end{equation}
Condition \eqref{nondeg} ensures that the map
\begin{equation}\label{eq:param-Leg}
\{ (y,v) \mid d_v\phi(y,v) = 0 \} \to \{ (y,\tau,\mu) \mid \tau =-\phi, \mu
= d_y \phi \}
\end{equation}
is a diffeomorphism.
The class $\Isc^m(X,L)$ of Legendre distributions are defined as a finite sum
of terms of the form
$$
x^{m+n/4-p/2} \int_{\RR^p} e^{i\phi(y,v)/x} a(x,y,v) dv,
$$
where $\phi$ locally parametrizes $L$ and $a$ is smooth, with compact
support in $v$. (We also allow the case that $\phi$ is a linear function of
$v$ and $a$ is Schwartz in $v$; these are called `extended Legendrian
distributions' in \cite{Hassell:Plane}.) It is not hard to see that the
above examples
are instances of this framework.

The {\it microsupport} of $u$ is the closed subset of $L$ corresponding,
under \eqref{eq:param-Leg}, to the set
$$
\{ (y,v) \mid d_v \phi = 0 \ \text{ and there is no neighbourhood of } (0,y,v)
\text{ in which } a \text{ is } O(x^\infty) \}.
$$
By a partition of unity we can always write $u$ as a finite sum of terms
each having microsupport as small as desired.

\

Such distributions turn up naturally in scattering theory. For example, let
$X$ be the radial compactifiction of $\RR^n$, and let $z$ be a linear
coordinate on $\RR^n$. Consider the function $e^{-iz \cdot k}$, $k \in
\RR^n$ which is a generalized eigenfunction of $\Lap$ with eigenvalue
$|k|^2$. In inverted polar coordinates $(x,\hat z)$, $x = |z|^{-1}$, $z =
\hat z/x$, this function takes the form $e^{i\hat z \cdot k/x}$, which is a
Legendre distribution of the first type. An example of a Legendre
distribution of the second type is a function of some subset of the $z$
variables: let $z = (z',z'')$ and let $V$ be a Schwartz function of
$z'$. Then
$$
V(z') = (2\pi)^{-k} \int e^{i\hat z' \cdot \zeta'/x} \hat V(\zeta') d\zeta'
$$
is a Legendre distribution of the second type. Such functions appear
naturally in the quantum $N$-body problem for example.

\

However, very frequently one comes across functions which are more
complicated, and associated to more than one Legendre submanifold. For
example, consider the kernel of the outgoing resolvent $R(\sigma + i0)$ of
the Laplacian on $\RR^3$. This is
$$
G(\sigma, z,z') = \frac1{4\pi} \frac{e^{i\lambda |z-z'|}}{|z-z'|}, \quad
\lambda = \sqrt{\sigma} > 0.
$$
Let us multiply by a function $\chi(|z-z'|)$ which is smooth, vanishes near
$|z-z'| = 0$ and is $\equiv 1$ for $|z-z'| > c$ to get rid of the interior
singularity. (The difference is the kernel of a pseudodifferential operator,
which is well understood.) Consider the resulting function, for fixed
$\lambda > 0$,
as a function on the radial compactification $\tilde X$ of $\RR^6$. Let $C
\subset
\pa \tilde X$ be the boundary of the diagonal $z=z'$. Away from $C$,
$|z-z'|$ is a
smooth function homogeneous of degree zero, so can be written $\phi(y)/x$,
where $y$ is a coordinate on $\pa \tilde X$ and $x = (|z|^2 +
|z'|^2)^{-1/2}$ is
the reciprocal of the Euclidean distance from the origin. Thus the kernel
is a Legendre distribution associated to a Legendre submanifold $L_1$ of
the first type in this region. On the other
hand, if we restrict to a region $|z-z'| < R$ near the diagonal then it
looks like a Legendre distribution associated with a Legendre submanifold
$L_2$ of the second type: a smooth function of
the three variables $z-z'$. Clearly it is simultaneously associated to both
of these Legendrian submanifolds. The geometry of these submanifolds is such
that $L_1$ is a manifold with boundary, which intersects $L_2$ cleanly at
$L_1 \cap L_2 = \pa L_1$.

There are (at least) two ways of looking at such distributions, and
consequently, two classes of distributions associated to $(L_2,L_1)$. The
first way is to define {\it intersecting Legendre distributions} as was
done in \cite{Hassell:Plane} (which is a routine generalization of the class of
intersecting Lagrangian distributions to the Legendre setting). One defines
a local parametrization of $(L_2,L_1)$ near a point $q \in L_1 \cap L_2$ to
be a function $\phi(y,v,s)$, where $v \in \RR^p$ and $s \in [0,\infty)$,
such that $\phi$ parametrizes $L_1$ in the sense analogous to
\eqref{eq:param-Leg} with both $v$ and $s$ taken as parameters,
while $\phi(y,v,0)$
parametrizes $L_2$. The nondegeneracy condition \eqref{nondeg} is replaced
by
\begin{equation}
d \left( \frac{\pa \phi}{\pa s}\right), d \left( \frac{\pa \phi}{\pa v_i} \right)
\text{ and } ds
\text{ are linearly independent. }
\label{nondeg2}\end{equation}
Then an expression of the form
\begin{equation}
x^{m+n/4-(p+1)/2} \int \int_0^\infty e^{i\phi(y,v,s)/x} a(x,y,v,s)\, ds\, dv,
\label{eq:int-Leg}\end{equation}
is a Legendre distribution of order $m$ associated to $L_1$ (e.g.\ when $a$ is
supported in $s \geq \epsilon > 0$) and of order $m+1/2$ at $L_2 \setminus
L_1$ (to see
this, multiply and divide by $d_s \phi/x$ and then integrate by parts in
$s$ to get a boundary term at $s=0$).
The class $\Isc^m(X,(L_2,L_1))$ is defined to be those functions $u = u_1 + u_2
+ u'$, where $u_1 \in \Isc^{m}(X,L_1)$, $u_2 \in \Isc^{m+1/2}(X,L_2)$ and $u'$
given by a finite sum of terms of the form \eqref{eq:int-Leg}.

The second way, when $L_2$ arises from an embedded submanifold $C \subset
\pa X$ as above, involves blowing up the submanifold $C$ and defining {\it
fibred Legendrian distributions}. Let
$$
Y = [X;C]
$$
be the manifold with codimension 2 corners obtained by real blowup of $C$.
$Y$ has two boundary hypersurfaces: one which is the lift of $\pa X$, which
we call the `main face', and a new hypersurface arising from the blowup of
$C$ which we call the `front face'. These will be abbreviated mf and ff,
respectively.
\begin{figure}\centering
\epsfig{file=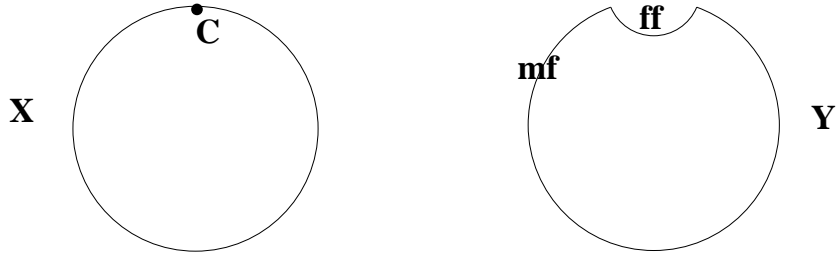,width=11cm,height=3.3cm}
\caption{Blowing up $C \subset \pa X$ to produce $Y$.}
\end{figure}
Then $\Isc^{m+1/2}(X,L_2)$ is identical with functions of the form
$x^{m+1/2-k/2+n/4} c$, where $c$ is smooth on $Y$ and vanishes to infinite
order at mf, and $k$ is the dimension of the fibers of the blow-down map,
i.e.\ it is the codimension of $C$ in $\pa X$.
Suppose $L_1$ is a Legendrian with boundary
meeting $L_2$ cleanly at $\pa L_1 = L_1 \cap L_2$. Let $(y',y'')$ be
coordinates on $\pa X$ which
define $C$ as $\{ x=0, y' = 0 \}$. If we assume that $L_1 \cap L_2$ has a
full rank projection to $C$, then this implies that near $q \in \pa L_1$,
there must be one of
the $y'$ coordinates, say $y'_k$, with a nonzero differential at $q$ (this
is shown in the following section). Then,
there is a parametrization of $L_1$ near $q$ of the form
$$
\phi(y,v)/x = y_k \tilde \phi (y'/y_k,y_k,y'',v)/x, \quad v \in \RR^p.
$$
We define the class of fibred Legendrians of order $(m,r)$ associated to
$L$ to be those functions $u = u_1 + u_2
+ u'$, where $u_1 \in \Isc^{m}(X,L_1)$ and supported away from $C$, $u_2 \in
\Isc^{r}(X,L_2)$ and $u'$
given by a finite sum of terms of the form
\begin{equation}
y_k^{r+\frac{n}{4}-\frac{k}{2}}
\left(\frac{x}{y_k}\right)^{m+\frac{n}{4} - \frac{p}{2}}
\int e^{iy_k \tilde
  \phi(\frac{y'}{y_k},y_k,y'',v)/x} a(y_k,\frac{y'}{y_k},
\frac{x}{y_k},y'',v) dv.
\label{eq:fibred-Leg}\end{equation}
Again $a$ is assumed to be smooth and compactly supported in all
variables. Note that this now means that $a$ is smooth in terms of the
differentiable structure on $Y$, rather than on $X$.

The purpose of this paper is to clarify the relation between these two
spaces $\Isc^m(X, (L_2,L_1))$ and $\Isf^{m,m+1/2}(Y,L)$ when $L_2$ is the
Legendre submanifold associated to $C$ as above.

\section{Invariant Description}
To describe the situation more invariantly,
let $X$ be a manifold with boundary, $x$ a boundary defining function,
and let $C$ be a closed embedded submanifold of $\pa X$. We denote the
interior of $X$ by $X^\circ$.
Let $p \in C$ and let $(x,y) = (x, y_1, \dots, y_{n-1})$ be coordinates
near $p$, where $x$ is a defining function for $\pa X$ and $C = \{ x=0, y'
= 0 \}$ near $p$; here $y' = (y_1, \dots, y_k)$.
Then $X$ is naturally equipped with its scattering cotangent bundle,
$\scT^*X$. One way to describe $\scT^* X$ is that its smooth sections
are spanned, over
$\Cinf(X)$, by one-forms of the form $d(\phi/x)$, $\phi\in\Cinf(X)$.
In particular, $\scT^*_{X^\circ}X$ is naturally identified with
$T^*X^\circ$. A local basis for $\scT^* X$ near $p$ is given by $dx/x^2 = -
d(1/x)$ and $dy_i/x$. Thus, a point $q \in \scT^* X$ may be written
\begin{equation}
q = \tau \frac{dx}{x^2} + \sum_i \mu_i \frac{dy_i}{x},
\label{eq:scT-coords}\end{equation}
and this gives local coordinates $(x,y,\tau,\mu)$ on $\scT^* X$, where
$(\tau, \mu)$ are linear coordinates on each fibre.
Moreover, $\scT^*_{\pa X}X$ is naturally
equipped with a contact structure
via a contact form $\chi$ induced by the symplectic structure of
$T^*X^\circ$, in a similar way to that in which a contact form is induced
on the cosphere bundle $S^*X=(T^*X\setminus 0)/\RR^+$.
In the local coordinates given by \eqref{eq:scT-coords}, this form
is equal to $\chi = d\tau + \mu \cdot dy$, as in the previous section.

There is a well-defined
`scattering conormal bundle' over $C$, denoted $\scN^*(C;X)$, which
is defined as the span of $d(\phi/x)$ (inside
$\scT^*_C X$) for all $\phi\in\Cinf(X)$
which vanish on $C$. In the local coordinates $(x,y',y'')$,
$\phi$ can be written
as $\phi=\sum_j y'_j a_j +x a_0$, so $\scN^*(C;X)$ is spanned by the
$dy'_j/x$, i.e.\ the dimension of each fibre is is the codimension of $C$
in $\pa X$.
It is easy to see that $\scN^*(C;X)$ is in fact a Legendre submanifold
of $\scT^*_{\pa X}X$.

Let
$Y=[X;C]$ be the blow-up of $X$ at $C$, and let $\beta:Y\to X$ be
the blow-down map. Let $\mf$ denote the lift of $\pa X$ to $Y$, and
let $\ff$ denote the front face of the blow-up, i.e.\ the lift of $C$.
Then $\ff$ has a natural fibration over $C$ given by the blow-down map:
$\beta|_{\ff}:\ff\to C$. As discussed in \cite{Hassell-Vasy:Spectral},
this defines a structure algebra $\Vsf(Y)$ of vector fields, and
more importantly for us, a corresponding replacement of the
standard cotangent bundle, namely the scattering fibred cotangent
bundle $\sfT^*Y$. Sections of $\sfT^*Y$ are spanned, over $\Cinf(Y)$,
by $d(\phi/x)$ where $\phi\in\Cinf(Y)$ is constant on the fibers of
the fibration. Such a setting is a natural generalization (to manifolds
with corners) of the
fibred cusp Lie algebra introduced by Mazzeo and Melrose on manifolds
with a fibred boundary \cite{Mazzeo-Melrose:Fibred}.
Then $\sfT^*_{\mf}Y$ has a natural contact form
which degenerates at the corner $\sfT^*_{\mf\cap\ff}Y$.
In this setting, $\sfT^*Y$ is just the pull-back of
$\scT^*X$ by the blow-down map $\beta$, and we denote the induced map
by $\tilde\beta:\sfT^*Y\to \scT^*X$. Thus, local coordinates on $\sfT^* Y$ near
the boundary consist of local coordinates on $Y$ together with the
functions $\tau$ and $\mu$ lifted from $\scT^* X$.
Moreover, the contact form on
$\sfT^*_{\mf}Y$ is just the
pull-back of $d\tau + \mu \cdot dy$ by $\tilde\beta$.
Given a Legendre submanifold $L\subset
\sfT^*_{\mf}Y$ which is transversal to the corner $\mf\cap\ff$ and
satisfies a compatibility condition with the fibration, the
class of Legendre distributions associated to $L$ was defined
in \cite{Hassell-Vasy:Spectral}; in the case of interest here, the
definition \eqref{eq:fibred-Leg} above suffices.
In particular, if $L$ has full rank projection to $\mf$, the compatibility
condition is automatically satisfied, and such a Legendre distribution is
simply of the form $e^{i\phi/x} a$ where $\phi\in\Cinf(Y)$ is constant
along the fibers of $\beta|_{\ff}$, and $a\in\Cinf(Y)$.

\section{Main Results}\label{sec:thm}
To describe our main results, we let $L_2$ be the scattering conormal
bundle $\scN^*(C;X)$ as above, and suppose that $L_1 \subset \scT^* X$
is a Legendre submanifold with boundary which intersects $L_2$ cleanly in
$\pa L_1 = L_2 \cap L_1$. In addition, we assume that $L_2 \cap L_1$ has a
full rank projection to $C$.
Thus, in local coordinates in $\scT^*_{\pa X}X$,
$L_2 = \{y'=0,\ \tau=0,\ \mu''=0 \}$. We shall prove in Lemma~\ref{contact}
that these assumptions imply that
$L_1$ also intersects $\scT^*_C X$ cleanly, with intersection $L_2\cap L_1$,
and in particular that at any $q\in L_1$, the pull-back of
$dy'_j$ to $L_1$ does not vanish for some $j$. Without loss of generality
we may assume that the restriction of $dy'_k$ to $L_1$ does not vanish at
$q$; thus, near $q$, we may assume that $y'_k \geq 0$ is a boundary
defining function for $L_1$.
This implies that under the blow-down map $\beta^*:\sfT^*_{\mf} Y\to
\scT^*_{\pa X}X$,
$L_1$ lifts to a Legendre submanifold $L$ which is transversal to
$\sfT^*_{\mf\cap \ff}Y$ and has a full rank projection to $C$.

Let $\Lt_1$ be an extension of $L_1$ across its boundary to an (open)
Legendre submanifold.
Our first result is a characterization of $\Isc^m(X,(L_2,L_1))$ in terms of
$I^m(X,\Lt_1)$.

\begin{thm}\label{thm:1}
Let $\alpha\in\Cinf(\RR)$, with $\alpha(t)$ identically $1$ for $t>1$,
identically $0$ for $t<0$. Let $U$ be a neighborhood of $q\in L_2\cap
L_1$ in $\scT^*_{\pa X}X$ such that $y'_k$ is a defining
function for $\pa L_1$ in $L_1\cap U$.
Intersecting Legendre distributions $u \in \Isc^m(X,(L_2,L_1))$
microsupported in $U$
may be represented as
\begin{equation}
\alpha\big(\frac{y_k}{x}\big) u_1 + u_0,
\label{eq:structure-int}
\end{equation}
where $u_1 \in \Isc^{m}(X,\Lt_1)$ and $u_0 \in
\Isc^{m+1/2}(X,L_2)$. Conversely, any such function is in
$\Isc^m(X,(L_2,L_1))$.
\end{thm}

Given this theorem we can rather easily understand the relation between
the spaces $\Isc^m(X,(L_2,L_1))$ and $\Isf^{m,m+1/2}(X,L)$.

\begin{thm}\label{thm:2} The space $\Isc^m(X,(L_2,L_1))$ is a proper subset of
$\Isf^{m,m+1/2}(X,L)$. In particular, $\Isf^{m,m+1/2}(X,L)$ is invariant
under multiplication by
$\Cinf(Y)$ functions, while $\Isc^m(X,(L_2,L_1))$ is not.
\end{thm}

{\it Remark. }One can show that the algebraic $\Cinf(Y)$ module generated
by the space $\Isc^m(X,(L_2,L_1))$ is dense in
$\Isf^{m,m+1/2}(X,L)$. Hence, one
could say that $\Isf^{m,m+1/2}(X,L)$ is generated over $\Cinf(Y)$ by
$\Isc^m(X,(L_2,L_1))$.

{\it Remark.} Guillemin and Uhlmann \cite{Guillemin-Uhlmann:Oscillatory}
and Joshi \cite{Joshi:Symbolic} have defined paired Lagrangian
distributions
of independent orders $(m,r)$ associated to a pair $(L_2, L_1)$ with the
geometry as described above.
Using the Fourier transform and a local identification of
neighbourhoods of $p \in \pa X$ with cones in $\RR^n$, one can define
paired Legendre distributions. In our setting these are defined
by
\begin{equation}
x^{m+n/4-(p+1)/2} \int \int_{-\infty}^\infty s_+^{r-(m+1/2)}
e^{i\phi(y,v,s)/x} a(x,y,v,s)\, ds\, dv,
\label{eq:int-Leg-indep}
\end{equation}
i.e.\ the distribution $s_+^{r-(m+1/2)}$ replaces the Heaviside step function
$H(s)=s_+^0$. It is not hard to show that Theorem~\ref{thm:2}
holds in this setting with distributions of the form
\eqref{eq:int-Leg-indep} forming a proper subset
of the class $\Isf^{m,r}(X,L)$.

{\it Remark. }It is very natural to assume that
$L_1\cap L_2$ is codimension one in $L_2$ and in $L_1$ since this appears
naturally in real principal type propagation. However, one can also
consider the case of higher codimension intersections. In the case of
Lagrangian manifolds, such intersections were studied by Guillemin
and Uhlmann \cite{Guillemin-Uhlmann:Oscillatory}. They essentially
define the associated class of distributions by blowing up the
intersection $L_2\cap L_1$. Note that if the intersection has codimension
one, the blow-up divides $L_1$ into two manifolds with boundary, each with
boundary $L_2\cap L_1$, hence resulting exactly in the setting discussed
above. In our case thus the most natural definition
of such intersecting Legendre distributions associated to an intersecting
pair $(L_2,L_1)$ with $L_2\cap L_1$ having codimension greater than one
is via the fibred scattering structure rather than directly by oscillatory
integrals. In particular, if $L_2=\scN^*(C;X)$, $L_1$ is the zero section
of $\scT^*_{\pa X}X$, the codimension of the intersection $L_2\cap L_1$
(in $L_2$ and in $L_1$) is given by
the dimension of the fibers of $\scN^*(C;X)\to C$, namely by the
codimension of $C$ in $\pa X$. The natural definition of the class
$\Isc^m(X,(L_2,L_1))$ of distributions
associated to the intersecting Legendre pair $(L_2,L_1)$ is
functions of the form
$x^{m+\dim X/4}\alpha(|y'|/x)f+g$, $f\in\Cinf(X)$, $g\in\Cinf([X;C])$ with
infinite order vanishing on $\mf$.

\section{An example}
Although the general case is hardly more complicated, for the sake of clarity
we first consider the case when $\dim X=2$, and $C\subset\pa X$ is a point
$x=0$, $y=0$, and $L_1$ is the zero section of $\scT^*X$ in $y\geq 0$.
The lift $L$ of $L_1$ to $\sfT^*_{\mf}Y$ is the zero section then.
Thus, fibred Legendrians in $\Isf^{m,m+1/2}(X,L)$
are $\Cinf$ functions on
$Y=[X;C]$, multiplied by powers of the boundary defining function,
i.e.\ functions of the form $x^{m+1/2} c$, $c\in\Cinf(Y)$.
Let us show that elements of $\Isc^m(X,(L_2,L_1))$ are functions of the
form $x^{m+1/2} c$, where $c\in\Cinf(Y)$ is of the form
\begin{equation}
c = \alpha\big(\frac{y}{x}\big) f(x,y) + g\big(x,\frac{y}{x}\big).
\label{eq:example-int}
\end{equation}
Here $f$ is smooth, $\alpha$ is as in Theorem~\ref{thm:1}, and $g(x,Z)$ is a
smooth function of $x$ taking values in Schwartz functions of $Z$. Thus,
$c$ is locally an arbitrary element of $\Cinf(Y)$ away from the corner, while
at the corner, its Taylor series is restricted so that it only has terms of
the form $y^j (x/y)^k$ with $j \geq k$.

An intersecting
Legendre distribution, $u\in \Isc^m(X,(L_2,L_1))$,
associated to these Legendrians is
one which can be written as $u=u_0+u_1+u'$,with
$u_0\in \Isc^{m+1/2}(X,L_2)$, $u_1\in\Isc^m(X,L_1)$, and $u'$
of the form
\begin{equation}\label{eq:int-Leg-8}
x^{m-1/2}
\int_0^\infty\int e^{i\zeta\cdot (y-\yb)/x} a(x,y,\zeta,\yb)\,d\zeta\,d\yb,
\end{equation}
where $a$ has compact support in $x,y,\yb$,
and is Schwartz in $\zeta$.
Directly from their definition, Legendre distributions in $\Isc^{m}(X,L_1)$ are
in fact of the form $x^{m+1/2}f(x, y)$, where $f$ is rapidly decreasing at
the boundary wherever $y < 0$. On the other hand,
it is easy to see that Legendrians in $\Isc^{m+1/2}
(X,L_2)$ are in $x^{m+1/2}\Cinf(Y)$,
vanishing to infinite order off the front face. Indeed, by definition,
such a distribution $u_0$ can be written as
$$
x^{m+1/2}\int e^{i\zeta y/x} a(x,y,\zeta)\,d\zeta,
$$
with $a$ compactly supported in $x,y$, rapidly decreasing in $\zeta$.
But this is just Fourier transform in $\zeta$, hence the result is
of the form $x^{m+1/2}b(x,y,y/x)$, with compact support in $x,y$, rapid decay
in $y/x$, which means that $u_0\in x^{m+1/2}\Cinf(Y)$ vanishing to
infinite order off the front face. Conversely,
if $u_0=x^{m+1/2}v$, $v\in\Cinf(Y)$ with infinite order vanishing on $\mf$,
then modulo a Schwartz function, $v$ is a $\Cinf$ function of $x$, $y$
and $y/x$ which is compactly supported in $x$ and $y$, and Schwartz in
$y/x$. Thus, we can
write $u_0$ (modulo $\dCinf(X)$)
as a Legendre function associated to $L_2$ by writing $v$
as the inverse Fourier transform of its Fourier transform in $y/x$:
$$
u_0(x,y)=(2\pi)^{-1}x^{m+1/2}\int e^{i\zeta y/x} \hat v(x,\zeta)\,d\zeta,
$$
modulo $\dCinf(X)$,
which is the definition of a Legendre distribution associated to $L_2$.

Hence the real question is whether pieces
as in \eqref{eq:int-Leg-8} are of the form \eqref{eq:example-int},
and conversely, whether such functions
supported near the corner can be written as in \eqref{eq:int-Leg-8} modulo
$\Isc^{m+1/2}(X,L_2)+\Isc^m(X,L_1)$.

Start with the former, i.e.\ consider an oscillatory integral of the
form \eqref{eq:int-Leg-8}. The phase function $\phi=\zeta\cdot
(y-\yb)$, is stationary with respect to $\zeta$ if $y=\yb$, and
it is stationary with respect to $\yb$ if $\zeta=0$.
Thus, $L'_1$ corresponds to the
set $\zeta=0$, $y=\yb$, $\yb\geq 0$, in the parameter space, while
$L_2$ corresponds to $\yb=0$, $y=0$. We may, modulo a Schwartz error,
assume that $a$ is independent of $\yb$. The reason is that we can expand
$a$ as a Taylor series around $y=\yb$:
$$
a(x,y,\yb,\zeta) = a(x,y,y,\zeta) + (y-\yb)a_1(x,y,\yb,\zeta).
$$
Since $(y-\yb) = \pa_\zeta \phi$, so by integrating by parts to get rid of
this factor on the second term we gain a power of $x$. Thus, iterating the
procedure and using asymptotic completeness we only have to deal with the
first term.

We start with the $\zeta$ integral, which simply takes a Fourier transform
of $a$
in $\zeta$ and evaluates the result at $-(y-\yb)/x$. We write $\hat a$
for the Fourier transform; so $\hat a$ is a Schwartz function in its
third variable (since it is the Fourier transform of a Schwartz function
in that variable). Indeed, we can arrange that $a$ is such that its
Fourier transform in $\zeta$ has compact support (and is smooth). Letting
$Z = y/x$ and $Z' = (y-\yb)/x$ and changing variables in the $\yb$ integral
we obtain
\begin{equation}\label{eq:int-Leg-16}
\begin{gathered}
x^{m-1/2} \int_0^\infty \hat a(x,y,-(y-\yb)/x) d\yb
= x^{m+1/2} \int_{-\infty}^Z \hat a(x,y,-Z') dZ' \\
= x^{m+1/2} \int_{-\infty}^Z \left[ d(x,y)\alpha'(Z') + b(x,y,Z')\right]dZ'
\end{gathered}
\end{equation}
where $d(x,y) = \int_{\RR} \hat a(x,y,-Z') dZ'$ is chosen so that the
integral of
$b$ in $Z'$ (over $\RR$) vanishes.
Thus, the integral of $b$ from $-\infty$ to $Z$ is a
Schwartz function of $Z$, so we have written this in the form
$$
d(x,y) \alpha(Z) + g(x,xZ,Z).
$$
Expanding $g$ as a Taylor series in the second variable we see that we have
expressed $u$ in the form \eqref{eq:example-int}.

Conversely, let $u$ be of the form \eqref{eq:example-int}. The $g$ term is
associated to $L_2$, so we only need deal with the first term. In fact, it
is clear that the class $\Isc^m(X,(L_2,L_1))$ is invariant under
multiplication by smooth functions, so we need only treat the $\alpha$
function.  To do this,
we write
\begin{equation}\begin{gathered}
\alpha(Z) = \int_{-\infty}^Z \alpha'(Z') dZ' \\
= x^{-1} \int_0^\infty \alpha'\big( \frac{y-\yb}{x} \big) d\yb  \\
= x^{-1} \int d\zeta \int_0^\infty e^{i(y-\yb)\zeta/x} \widehat{\alpha'}(\zeta)
d\yb d\zeta
\end{gathered}\end{equation}
which is of the right form since $\alpha'$, and therefore also its Fourier
transform, are Schwartz.

\section{The general case}
Let $X$, $C$, $Y$, $L_2$, $L_1$, $\Lt_1$, and $L$ be as in
section~\ref{sec:thm},
and let $q \in \pa L_1$. Let $\Lt_1$ be a Legendrian extension of $L_1$ to
a submanifold without boundary across $L_2$.
In local coordinates, $C = \{ x=0, y'=0 \}$ and $L_2 = \{ x=0, y'=0,
\mu''=0 \}$. First we prove a statement asserted just before
the main results from section~\ref{sec:thm}.

\begin{lemma}\label{contact} There is a $y'$ coordinate, which may be taken
to be $y_k$
without loss of generality, whose differential restricted to $L_1$ does not
vanish at $q$.
\end{lemma}

\begin{proof} Let $q \in L_2 \cap L_1$, and let
$W=T_q (L_2\cap L_1)\subset Z=T_q\scT^*_{\pa X}X$. Moreover,
let $V_j=T_q L_j\subset Z$. The fact that a subspace $V$
of $Z$ is Legendre means that both $\chi$ and $d\chi$
vanish on it identically, i.e.\ $\chi(v)=0$, $d\chi(v,v')=0$ for all $v,v'
\in V$, and $V$ is maximal, i.e.\ $(\dim Z-1)/2=n-1$ dimensional,
$n=\dim X$, with this
property. Note that $d\chi$ is non-degenerate on $\Ker\chi\subset
Z$, i.e.\ it is a symplectic form on this vector space.

Now, both $\chi$ and $d\chi$ vanish on $W$ since $W\subset V_2$.
Let $W'$ denote the subspace of $\Ker\chi$ which annihilates $W$,
i.e.\ $W'=\{w'\in Z:
\ \chi(w')=0,\ d\chi(w',w)=0\ \text{for all}\ w\in W\}$. Then any Legendre
subspace $V\supset W$ of $Z$ satisfies $V\subset W'$
since $V\subset\Ker\chi$ and $d\chi(v,w)=0$ for all $w\in W\subset V$.
Note that $W$ has codimension $2$ in $W'$ (since $W$ has codimension $1$
in the Legendre subspace $V_2$). Thus, $W'/W$ is a 2-dimensional vector
space, and $d\chi$ descends to a symplectic form on it. The image $V'$ of a
Legendre subspace $V\supset W$ in $W'/W$ is Lagrangian with respect to
this form. There is a one-dimensional family of such Lagrangian
subspaces; the image of $V_2$ is one of them. Indeed,
given any non-zero element $u$ of $W'/W$, there is a unique Lagrangian subspace
of $W'/W$ which includes $u$, namely the span of $u$. This then determines
a unique Legendre subspace $V$ of $W'$ with $W\subset V$.

We claim that $W'$ is not a subspace of $T_q\scT^*_C X$. Indeed, suppose
otherwise, i.e.\ that $W'\subset T_q\scT^*_C X$. The hypothesis on the
full rank projection of $L_2\cap L_1$ means that $dy''_j$ are independent
on $W$. The corresponding Hamilton vectors $\pa_{\mu''_j}$ under $d\chi$ in
$\Ker\chi$ are tangent to $\scT^*_C X$, hence in $T_q\scT^*_C X
\cap\Ker\chi$. Thus, they span a $(\dim C)$-dimensional subspace $T$ of this
space. If $f$ is a nonzero linear combination of the functions $y''_j$,
then $df$ does not vanish on $W$, which implies that $d\chi(H_f, \cdot)$
does not vanish on $W$. Hence $H_f \notin W'$, which means that $T$ and
$W'$ have trivial intersection. But this is a contradiction:
by dimension counting, the codimension of $W'$ in $T_q\scT^*_C X
\cap\Ker\chi$ is $\dim C-1$. Hence $W'$ is not a subspace of $T_q\scT^*_C X$.

Thus $dy'_j$ cannot all vanish identically on $W'$.
Since they all vanish identically on $V_2$, this is the only Legendre
subspace with this property. By the clean intersection assumption,
$V_2\cap V_1=W$, i.e.\ $V_2$ and $V_1$ are not the same. Hence the
$dy'_j$ do not all vanish on $V_1$, i.e.\ the pull-back of $dy'_j$
to $L_1$ at $q$ is non-zero for some $j$. By relabelling the coordinates,
we may assume that $dy'_k$ is non-zero.
\end{proof}

Since $\pa L_1$ has a full rank projection to $C$,
the span of the pull-back of the differentials $dy'_j$ to $L_1$ at $q$ is
exactly one-dimensional. By a linear change of the $y$ coordinates we may
assume that $dy'_1,\ldots,dy'_{k-1}$ pull back to $0$ at $q$. Then
standard contact arguments show that
$y'_k$, $y''$ and $v = (\mu'_1,\ldots,\mu'_{k-1})$ give local coordinates on
$\Lt_1$ near $q$, and $L_2\cap L_1$ is defined by $y_k=0$ in these
coordinates. By switching the sign of $y_k$ if necessary, we may also assume
that microlocally $L_1$ lies in $y_k\geq 0$.
Expressing $\yt =(y_1,\ldots,y_{k-1})$ as
$\yt =\Yt(y_k,y'',v)$, $\tau$ as $\tau=T'(y_k,y'',v)$ on $\Lt_1$,
it follows that a local nondegenerate
parameterization of $\Lt_1$ is given by
\begin{equation}
\phi(y,v)=-T'+v \cdot(\yt - \Yt), \quad v \in \RR^{k-1},
\label{eq:phi}\end{equation}
while a local parametrization of $L$ is given by the same function for
$y_k/x > C$. A local nondegenerate parametrization of $(L_2,L_1)$ is given by
\begin{equation}
\psi(y,v,\zeta,\yb)=-T'+v \cdot(\yt - \Yt)+\zeta(y_k-\yb), \quad \yb \geq 0.
\label{eq:psi}\end{equation}
Note that $\yt$ when $y_k=0$, so in fact $\Yt=y_k Y$, and similarly
$T'=y_k T$ (since $\tau=0$ on $L_2$).

\

{\it Proof of Theorem~\ref{thm:1}. }
Write $u \in \Isc^m(X,(L_2,L_1))$ in terms of the phase function $\psi$
from \eqref{eq:psi}.
We may then run the argument in the example of the previous
section in the variables $(x,y_k,\yb,\zeta)$ to show that $u$ can be
written in the form \eqref{eq:structure-int}. To prove the converse, write
$u_1$ with respect to the phase \eqref{eq:phi} and express $\alpha$ as in
\eqref{eq:int-Leg-16} to obtain an expression involving the phase function
$\psi$ from \eqref{eq:psi}, which is manifestly an element of
$\Isc^m(X,(L_2,L_1))$.

\

{\it Proof of Theorem~\ref{thm:2}. }
First we show that $\Isc^m(X,(L_2,L_1))$ is contained in
$\Isf^{m,m+1/2}(X,L)$.
Theorem~\ref{thm:1} to write $u \in \Isc^m(X,(L_2,L_1))$
in the form
\eqref{eq:structure-int}. In terms of this representation, $u_0$ is
$x^{m+n/4-(k-1)/2}$
times an element of $\Cinf(Y)$ which vanishes to all orders at the main
face, so this is certainly an element of $\Isf^{m,m+1/2}(X,L)$ (in fact,
$\Isf^{r,m+1/2}(X,L)$ for any $r$). On the other hand, $u_1$ is in $\Isc^m(X,
\Lt_1)$ so can be written with respect to the phase function $\phi$ from
\eqref{eq:phi}. But this is also a phase function for $L$, and
multiplication by a suitable $\alpha$ means it is now supported in $y'_k/x
> C$, so $\alpha \cdot u_1$ is also in $\Isf^{m,m+1/2}(X,L)$.

To show that the inclusion is proper, we write $u$ in terms of a reduced
symbol and show that its Taylor series is restricted at $x/y'_k = y'_k =
0$, that is, at the intersection of the front face and the main face on
$Y$. Thus, we can write $u$ in the form \eqref{eq:fibred-Leg} where $a$
only
depends on $y_k, x/y'_k, y''$ and $v$; also note that when $r=m+1/2$ and
$p=k-1$, then the power of $y'_k$ outside the integral vanishes. By the
symbol calculus of \cite{Hassell-Vasy:Resolvent}, the reduced symbol is then
determined to all orders in Taylor series at $\mf \cap \ff$ by $u$.
Using the description given by
Theorem~\ref{thm:1}, we see that the $u_0$ term has
trivial Taylor series at $x/y'_k = y'_k = 0$. The $\alpha \cdot u_1$ term
has the property that the symbol for $u_1$ is smooth in the variables $x$
and $y'_k$, so the Taylor series of $a$ as a function of $x/y'_k$ and
$y'_k$ has the property
\begin{equation}
\text{ The coefficient of } (x/y'_k)^j y_k^l \text{ vanishes whenever } l
< j.
\label{eq:Taylor}\end{equation}
Thus, the sum of the two terms $u_0$ and $\alpha \cdot u_1$ has property
\eqref{eq:Taylor}. It is
clear that this property is not invariant under multiplication by
smooth functions of $y'_k$ and $x/y'_k$. However, the space
$\Isf^{m,m+1/2}(X,L)$ is by its definition invariant under $\Cinf(Y)$. So
$\Isc^m(X,(L_2,L_1))$ is a strictly smaller space than
$\Isf^{m,m+1/2}(X,L)$.
This completes the proof of the Theorem.

\

{\it Remark. }The absence of terms \eqref{eq:Taylor} in the symbol of
$\Isc^m(X,(L_2,L_1))$ is reflected in the symbol calculus for $u \in
\Isc^m(X,(L_2,L_1))$. That is, if the symbol $\sigma^m_{L_1}(u)$ of $u$
on $L_1$ vanishes, then $u \in \Isc^{m-1}(X,(L_2,L_1)) +
\Isc^{m+1/2}(X,L_2)$ (see \cite{Melrose-Uhlmann:Intersection}, equation
(5.2)). On the other hand, if the symbol of $u \in
\Isf^{m,m+1/2}(X,L)$ at $L$ vanishes, then $u \in \Isf^{m-1,m+1/2}(X,L)$
(as opposed to $\Isf^{m-1,m-1/2}(X,L) + \Isf^{\infty,m+1/2}(X,L)$). This
better vanishing property of intersecting Legendre distributions makes them
more useful for understanding principal type propagation.

\bibliographystyle{plain}
\bibliography{sm}

\end{document}